\documentclass[12pt]{amsart}

\usepackage{graphicx}
\usepackage{amssymb}
\usepackage{amscd}

\newtheorem{thm}{Theorem}[section]
\newtheorem{prop}[thm]{Proposition}
\newtheorem{lem}[thm]{Lemma}
\newtheorem{cor}[thm]{Corollary}

\newtheorem*{thm*}{Theorem}
\theoremstyle{definition}

\newtheorem*{rmk}{Remark}

\newcommand{\bq}{\mathbb{Q}}
\newcommand{\bz}{\mathbb{Z}}
\newcommand{\br}{\mathbb{R}}
\newcommand{\bc}{\mathbb{C}}

\newcommand{\lan}{\left\langle}
\newcommand{\ran}{\right\rangle}

\newcommand{\rk}{\operatorname{rk}}
\newcommand{\fix}{\operatorname{Fix}}

\renewcommand{\hom}{\operatorname{Hom}}

\newcommand{\trace}{\operatorname{Trace}}
\newcommand{\charac}{\operatorname{char}}

\begin{document}

\title{Symmetry groups of non-simply connected four-manifolds }
\author{Michael McCooey}\date{\today}

\begin{abstract} Let $M$ be a closed, connected, orientable topological four-manifold with $H_1(M)$  nontrivial and free abelian, $b_2(M)\ne 0, 2$, and $\chi(M)\ne 0$. We show that if $G$ is a finite group of 2-rank $\le 1$  which admits a homologically trivial, locally linear, effective action on $M$, then $G$ must be cyclic.  

With additional assumptions to ensure orientability of some components of the singular set (e.g. if $G$ acts by symplectic symmetries, or preserving a spin structure), we also rule out $C_2 \times C_2$ actions. 

The proofs use equivariant cohomology, localization, and a careful study of the first cohomology groups of the (potential) singular set. 
\end{abstract}\maketitle
\section{Introduction}\thispagestyle{empty}

This paper can be viewed as a sequel to \cite{MM2}, where, following a conjecture of Edmonds~\cite{Actions}, we 
showed that if $M$ is a simply-connected four-manifold with $b_2(M)\ge 3$, and $G$ is a finite or compact Lie group which acts effectively, locally linearly, and homologically trivially on $M$, then $G$ must be isomorphic to a subgroup of $S^1\times S^1$. 

Here we consider the more general situation in which $H_1(M;\bz)$ is free abelian of arbitrary rank.  We prove:

\begin{thm}[Main Theorem]\label{MainTheorem}
Let $M$ be a closed, connected, orientable topological four-manifold with $H_1(M)$  nontrivial and free abelian, $b_2(M)\ne 0, 2$, and $\chi(M)\ne 0$. If $G$ is a finite group of 2-rank $\le 1$  which admits a homologically trivial, locally linear, effective action on $M$, then $G$ is cyclic. 
\end{thm}

The assumptions that $\chi(M)\ne 0$ and $b_2(M)\ne 2$ are necessary, as familiar examples of actions on $T^4$ and the product of $S^2$ with any closed, oriented surface make clear. We conjecture that the assumption that $G$ has $2$-rank 1 is not necessary, but the possibility of nonorientable surfaces in the singular set makes the problem more difficult. 
Although we do not explicitly indicate so each time, many arguments extend to the case where torsion in $H_1(M)$ is relatively prime to the orders of the groups involved. 

The analysis of ~\cite{MM2} was based on a comparison of the Borel equivariant cohomology of $M$ with that of its singular set $\Sigma$, and an important ingredient in understanding $\Sigma$ was Edmonds's observation~\cite{Aspects} that in the simply-connected case, the fixed-point set of any cyclic group action consists only of isolated points and spheres, with no surfaces of higher genus. Using local considerations, it becomes possible to assemble these fixed sets to gain a rather explicit description  of $\Sigma$ for a (potential) action of a larger finite group, and eventually rule out the nonabelian ones.
The arguments were homological in nature, and hence extend to the case where $H_1(M;\bz)=0$.

Our current situation differs in several important ways. First, algebraic considerations now allow (indeed, force) fixed-point sets to contain surfaces of higher genus, and a purely combinatorial assembly of $\Sigma$ from its components is no longer feasible. Moreover, when all cohomology of $M$ was concentrated in even degrees, nearly complete information about the cohomology of $\Sigma$ could ultimately be extracted from the map $H^2(M)\to H^2(\Sigma)$. In contrast, for a cyclic group $C_p$, the restriction $H^1(M)\to H^1(M^{C_p})$ has half rank at best, so it now requires more work to understand and exploit the induced action of a group $G$ on the cohomology of its singular set. Finally, in ~\cite{MM2}, gaps in odd degrees frequently led the spectral sequences involved in computing $H_G^*(M)$ to collapse for formal reasons.  They require more careful attention here, and in some cases, our understanding of the differentials remains incomplete. On the positive side, $H^1(\Sigma)$, to the extent that it can be detected, ultimately imposes considerable rigidity on a group action. 

The structure of the paper is as follows: In section~\ref{General}, we collect and provide references for our main technical tools: the Lefschetz Fixed-point Theorem, Borel equivariant cohomology, and the  Localization Theorem of Borel, tom Dieck, Hsiang, and Quillen. In section~\ref{Collapsing}, we consider the collapsing of the Borel Spectral Sequence for various groups and derive some immediate consequences, including generalizations of a few useful results of~\cite{Aspects}. In section~\ref{SigmaAction}, we  discuss, modify, and generalize some (slightly flawed) arguments of~\cite{Assadi} to analyze the action of a rank two abelian group on the fixed set of a cyclic subgroup, and ultimately rule out actions by most groups of rank two. We then discuss some of the subtleties of the case $p=2$.  Finally, in section~\ref{nonabelian}, we consider actions of metacyclic and quaternion groups, rule them out (using an application of localization which might be of independent interest), and gather the pieces to prove the main theorem.

Our standing assumptions throughout the paper are that $M$ is a closed, connected, orientable topological four-manifold, and $G$ is a finite group acting effectively, but with trivial induced action on $H_*(M;\bz)$,  and locally linearly (``HTLL") on $M$.  $C_n$ denotes $\bz/n\bz$ regarded as a (transformation) group. Cohomology should be assumed to be singular cohomology with integer coefficients in the absence of indications to the contrary.

\section{General considerations and background}\label{General}

Let $C_n=\lan g\ran$. Local linearity implies that $M^g$ is a locally flat submanifold of $M$, and if the action preserves orientation on $M$ (as we shall always assume), the fixed-point set must have even codimension, and hence consist of a union of points and surfaces. When $n$ is odd, each surface component must be orientable (see ~\cite{Bredon}). For locally linear actions, we also have a strong form of the Lefschetz fixed-point theorem: $\chi(M^g)=\sum (-1)^i\trace(g^*|_{H_i(M; \bq)})$. If $g$ acts trivially on homology, it follows that $\chi(M^g)=\chi(M)$.

Recall that the Borel construction $M_G =M\times_G EG$ defines a fibration
$M\to M_G\to BG$. We refer to the Leray-Serre spectral sequence of this fibration
$E_2^{ij}(M)=H^i(G; H^j(M))\Rightarrow H^*(M_G)$ as the \emph{Borel Spectral Sequence} (BSS), and  the groups 
$H_G^*(M):=H^*(M_G)$ as \emph{(Borel) equivariant cohomology groups.} For more background, see any of \cite{Borel}, \cite{tomDieck}, \cite{Bredon}, or \cite{Hsiang}.
	
	Several facts about these groups and the spectral sequence are of particular importance for us:
	\begin{enumerate}
	\item{The BSS is equipped with a  well-behaved  $H^*(G)$-algebra structure. In particular, if $H^*(M)$ is torsion-free and $G$ acts trivially on it, then 
	$E_2\cong H^*(M)\otimes H^*(G)$.}
	\item{Let $\Sigma :=\{ x\in M\ |\ G_x\ne \{e\} \ \}$ denote  the \emph{singular set} of the action. Then restriction induces an isomorphism 
	$H_G^*(M)\to H_G^*(\Sigma)$ in dimensions $*>4$ (or more generally, greater than the dimension of the manifold under consideration).}
	\item{The functor $H_G$ is natural with respect to maps of groups and $G$-spaces.} 
	\end{enumerate}

Finally, recall the Localization Theorem (cf. ~\cite{tomDieck, AlldayPuppeBook}): Let a finite group $G$ act on a compact space $X$, let $S$ be a multiplicatively closed, central subset of $H^*(G)$, and let $\Sigma_S=\{ x\in X\ |\ S\cap\ker(i^*:H^*(G)\to H^*(G_x)=0\}$. Then inclusion induces an $S^{-1}H^*(G)$-algebra isomorphism $S^{-1}H^*_{G}(X)\to S^{-1}H^*_G(\Sigma_S)$. Note that if $S$ includes an element  $a$ of degree $d$, then multiplication by $a$ is an isomorphism between   $S^{-1}H^*_G(X)$ and  $S^{-1}H^{*+d}_G(X)$,  so after localization, meaningful grading distinctions survive (at most)  only modulo $d$.

\section{Collapsing of the spectral sequence}\label{Collapsing}

Poincar\'e duality and the presence of a nonempty fixed-point set together impose strong restrictions on the differentials in the spectral sequence:
\begin{prop}\label{CollapsingLemma1}
Let $R=\bz$ or $\bz_p$. 
Suppose a finite group $G$ acts on a closed four-manifold $M$, with $H_*(M; R)$  $R$-torsion-free, trivially on $H^*(M;R)$, and  with at least one fixed point.  Then the differential $d_2$ in the BSS vanishes in all of row 4, all of row 3, all of row 1, and on those classes
in row 2 which are products of one-dimensional classes. Indeed, the only potentially nonzero differential $d_r$  (for $r\ge 2$) is $d_2^{i,2}$, and its values are determined by $d_2^{0,2}: H^0(G; H^2(M))\to H^2(G; H^1(M))$. If $d_2=0$, then the spectral sequence collapses. 

In particular, if $G$ is cyclic of prime order, then this spectral sequence collapses. 

\end{prop}

\begin{proof}
Let $x\in M^G$. The restriction homomorphism $j^*:H^*(M)\to H^*(M \setminus \{x\})$ is zero in dimension four, and an isomorphism in other dimensions. It follows that the corresponding map of spectral sequences $E_2^{i,j}(M)\to E_2^{i,j}(M\setminus\{x\})$ is trivial when $j=4$, and an isomorphism otherwise. Factoring through these maps, using naturality of the differentials, shows that $d_2^{i, 4}=0$ for all $i$.

Similarly, factoring through the map of spectral sequences $E_2^{i,j}(M, \{x\})\to E_2^{i,j}(M)$ shows that $d_2^{i, 1}=0$ for all $i$. 

If $H_*(M)$ is $R$-torsion-free, then the universal coefficient theorem and Poincar\'e duality together yield a nonsingular intersection pairing $H^*(M)\otimes H^{4-*}(M)\to H^4(M)$. So for any generator $u\in E_2^{0,3}=H^0(G; 
H^3(M))$, there is $v\in E_2^{0,1}=H^0(G;H^1(M))$ such that $uv$ generates $E_2^{0,4}=H^0(G;H^4(M))$. But $d_2(uv)=0 = ud_2(v)-d_2(u)v$. Since $d_2(v)=0$, $d_2(u)v=0$. But the $E_2$ term of the spectral sequence is a free $H^*(G)$-module, so $d_2(u)=0$, as well. 

Finally, if $G=C_p$, Adam Sikora~\cite[3.2(i), 3.11, 3.13]{SikoraTorusZp} has shown\footnote{Sikora's  result is stated for the Leray spectral sequence of the map $X_G\to B_G$, which is is defined in terms of sheaf cohomology. But according to Bredon~\cite[III.1.1]{BredonSheafTheory},whenever $X$ is a CW complex, $\Phi$ is a paracompactifying family of closed sets, and $\mathcal{A}$ is a sheaf of local coefficients,  there is a natural isomorphism $H^*_{\Phi}(X;\mathcal{A})\to {_{\Delta}}H^*_{\Phi}(X;\mathcal{A})$ between sheaf and ordinary singular cohomology groups. This isomorphism covers all cases of interest to us here. In particular, since the Leray and the Serre spectral sequences have isomorphic $E_2$-terms and abutments, and the Leray spectral sequence collapses, it follows for dimension reasons that the Serre spectral sequence must collapse, allowing us to bypass any technical verification that the two spectral sequences are themselves isomorphic. 

}  that the terms $E_r^{ij}$ of the spectral sequence satisfy a form of Poincar\'e duality for each $r\ge 2$ which implies that $\rk{E_3^{2,1}}=\rk{E_3^{2,3}}$, and hence that $d_2^{0,2}=0$.

Vanishing of $d_3$ on row 3 follows from the fact that generators of $H^3(M)$ are dual to those of $H^1(M)$, and all remaining differentials are easily accounted for by factoring through $E(M-\{x\})$ and $E(M, \{x\})$, as appropriate. 
\end{proof}

\begin{cor}\label{CollapsingLemma2}
Let $p$ be prime. In the above situation (in particular, with $M^G\ne \emptyset$), if $G=C_p\times C_p$, the BSS collapses with integral coefficients.
\end{cor}
\begin{proof}

Recall (cf. ~\cite{Lewis}) that 
$$H^*( C_2\times C_2; \bz)\cong\frac{\bz[\alpha_2, \beta_2]\otimes P[\mu_3]}
{\lan 2\alpha = 2\beta =2\mu =0, \mu^2=\alpha\beta^2+\alpha^2\beta\ran},$$
while for $p$ odd,
$$H^*( C_p\times C_p; \bz)\cong\frac{\bz[\alpha_2, \beta_2]\otimes \bigwedge[\mu_3]}
{\lan p\alpha = p\beta =p\mu =0\ran}.$$

In both cases, the elements $\alpha_2$ and $\beta_2$ arise as Bocksteins of elements of $H^1( C_p;\bz_p)=\hom( C_p, \bz_p)$, so it is easy to check that every $x\in H^2( C_p\times C_p;\bz)$ is detected by restriction to some cyclic subgroup $ C_p<G$. 

Hence, if $d_2^{0, 1}(x b)=ya\ne 0$ for some $x\in H^0( C_p\times C_p)$, $b\in H^2(M)$, $y\in H^2( C_p\times C_p)$, and $a\in H^1(M)$, then there is a restriction $r^*$ to some cyclic subgroup $ C_p< C_p\times C_p$ so that  $r^*(y)\ne 0$. But  $d_2$ commutes with $r^*$, so for the cyclic group action, $d_2(xb)\ne 0$, contradicting Proposition~\ref{CollapsingLemma1}. 
\end{proof}

It follows from the Universal Coefficient Theorem that all torsion in $H^*(M; \bz$ is determined by $H_1(M)$.
Henceforth we assume that $H_1(M)$ is torsion-free, and note some consequences of Lemma ~\ref{CollapsingLemma1}

\begin{cor}
Suppose $M$ is a four-manifold with $H_1(M)$ torsion-free, and suppose $ C_p$, where $p$ is prime, acts HTLL on $M$. 
If $\chi(M)\ne 0$, then the $\bz_p$-Betti numbers satisfy $b_1(M^{C_p})=2b_1(M)$, and $b_0(M^{C_p})+ b_2(M^{C_p})=2+b_2(M)$.
\end{cor}
\begin{proof}
In this situation, it follows from the Lefschetz fixed-point theorem that $\chi(M^{C_p})=\chi(M)\ne 0$. In particular, $M^{C_p}\ne \emptyset$.

Since  $H_G^*(M)\cong H_G^*(M^{C_p})$ in high degrees, it follows easily for odd $p$ that
$b_1(M^{C_p})=2b_1(M)$, $b_0(M^{C_p})+b_2(M^{C_p})=2+b_2(M)$. The case $p=2$ is slightly complicated by the possibility of 2-torsion in $H^2(M^{C_p}; \bz)$, but follows from~\cite[VII.3.1]{Bredon}, since the total Betti numbers of $M$ and $M^{C_p}$ agree. 
\end{proof}

Recall that a \emph{pseudofree} action of a finite group is one in which the singular set consists only of isolated points. We note in passing: 

\begin{cor}Suppose $M$ is a four-manifold with $H_1(M)$ torsion-free, and suppose $b_1(M)>0$. Then $M$ admits no pseudofree, homologically trivial group actions if $\chi(M)\ne 0$. If $\chi(M)=0$, then the only possible pseudofree actions are actually free actions.
\end{cor}

The following proposition is due to Edmonds\cite[2.5]{Aspects} in the simply-connected case. Generalization to the case at hand presents no new difficulties:

\begin{prop} \label{NontrivialRestriction} Let $M$ be a four-manifold with $H_1(M; \bz)$ torsion-free, with an HTLL action by $G= C_p$ ($p$ prime).  
If $F'$ is any proper subset of the fixed-point set $F$, then the restriction map $H^2(M; \bz_p) \to H^2(F' ; \bz_p)$ is surjective. 
\end{prop}
\begin{proof}
Without loss of generality, we assume that $F$ contains at least two points (and hence that the BSS collapses).  Let $y\in F-F'$, and $x\in F$. As in \cite{Aspects}, the map of spectral sequences $E(M-y, x)\to E(F-y, x)$ converges to an isomorphism in degrees $>4$, and in particular, in degree 6. The hypotheses on $H^1(M)$ imply that
$H^5(G; H^1(M; \bz))$ and $H^3(G; H^3(M;\bz))$ vanish, so $H^6(M_G-y_G; x_G)\cong H^4(G; H^2(M-y, x))\cong H^2(M; \bz_p)$. Meanwhile, the Leray-Hirsch theorem implies that the restriction $H^6(F_G-y_G, x_G)\to E^{4,2}_{\infty}(F-y, x)$ is onto. But $ E^{4,2}_{\infty}(F-y, x)\cong H^4(G; H^2(F-y, x; \bz))\cong H^2(F-y, \bz_p)$. 

\end{proof}

\begin{rmk}
This result actually generalizes to the case when $C_p$ acts nontrivially, but without cyclotomic-type summands, on $H^*(M)$.
One appeals to a stronger version of Sikora's result to show the spectral sequence still collapses, then shows that $H^5(G; H^1(M))$ and $H^3(G; H^3(M))$ make no contribution to $H^6_G(M)$. However, simple constructions show the result to be false when cyclotomic actions on $H^1(M)$ are permitted. 
\end{rmk}

 \begin{cor}[\cite{Aspects}] \label{InvariantSurfaces}In the above situation:
 	\begin{enumerate}
	\item{If $F$ is not purely $2$-dimensional, then the surface components of $F$ represent independent elements of $H_2(M; \bz_p)$.}
	\item{If $F$ is purely two-dimensional, with $k$ components, then the surfaces in $F$ span a subspace of $H_2(M;\bz_p) $ of dimension at least $k-1$, with any $k-1$ components representing independent elements. }
	
	\item{If $M^{C_p}$ has more than two components, then the surface components all represent different elements of $H_2(M)$, and hence can not be interchanged by a homologically trivial symmetry.}
	\end{enumerate}
\end{cor}

\section{The action of $G$ on $\Sigma$ and its consequences for rank two groups.}\label{SigmaAction}

Whenever $ C_p\triangleleft G$, there is an induced action of $G/ C_p$ on $M^{C_p}$. In this section we consider the extent to which this action can can be nontrivial.

Consider the following basic example. Let $G= C_2\times C_2=\lan h,k\ran$ act linearly on $S^2$. 
Then $\fix(h)$ consists of two points, and the $k$-action interchanges them. (Examples of actions on four-manifolds can be constructed in a similar vein, e.g the action on $S^2$ of any subgroup of $SO(3)$ defines obvious actions on $S^2$ bundles over surfaces, and  $ C_3\times C_3$ acts linearly on $\bc P^2$, with permutation actions on the fixed sets of cyclic subgroups. )

A slightly flawed argument of~\cite[Lemma 5.2]{Assadi} asserts that in such situations, $G/ C_p$ must act trivially on $H^*(M^{ C_p})$. The idea boils down to the following: In large dimensions, $H^*_{ C_p}(M)\cong H^*_{ C_p}(M^{ C_p})$. The action on the terms of the spectral sequence for the $ C_p$ action on $M$ is trivial;  so the action on $H^*(M^{ C_p})$ is claimed to be trivial as well.  

However,  $E(M)$ converges to the graded module associated to $H_{ C_p}(M)$ via the filtration of $H^*(G)$ by inclusion of skeleta of $BG$ -- not to the equivariant cohomology itself. Even with field coefficients, the action can be trivial on the spectral sequence, but nontrivial on the module. But analysis of the $G$-action on $E(M)$ can still be revealing. 

In the case of  $S^2$, let $H$ denote the subgroup $\lan h\ran$. The fact that the spectral sequence has only two rows yields the short exact sequence:
$$0\to E^{4,0}_{\infty}\to H^4_H(S^2)\to E^{2,2}_{\infty}\to 0.$$
The collapsing of the spectral sequence for the $H$-action, together with the isomorphism induced by the inclusion of the fixed-point set, means that the sequence becomes:
$$0\to H^4( C_2; H^0(S^2))\to H^4_H((S^2)^h)\to H^2( C_2;H^2(S^2))\to 0,$$
and finally, since the BSS for the fixed set has only one row, $H^4_H((S^2)^h)\cong H^4( C_2)\otimes H^0((S^2)^h))$, so 
we get
$$0\to\bz_2\otimes H^0(S^2)\to \bz_2\otimes H^0((S^2)^h)\to \bz_2\otimes H^2(S^2)\to 0.$$

As $h$ is (obviously) central in $G$,  the action of $G/H$ on $(S^2)^h$ is $H$-equivariant. 
It follows that $G/H$ acts on the spectral sequence as well, so the above is a short exact sequence of $\bz_2[G/H]$-modules, with trivial modules on the ends. 

Maschke's Theorem (see ~\cite{CurtisReiner}) provides one sufficient condition for triviality on the ends of such a sequence to imply triviality in the middle: an $\mathbb{F}[G/H]$-module $M$ 
is completely reducible whenever $|G/H|$ is relatively prime to $\charac(\mathbb{F})$.
 (Of course, in this case, the theorem does not apply.)  Here, the image of the map $\bz_2\otimes H^0(S^2)\to \bz_2\otimes H^0((S^2)^h)$ is the subgroup generated by $[x_0]+[x_1]$. The actions on both the subgroup and the resulting quotient are trivial.

This argument also breaks down for other reasons when $C_p$ is not central in $G$, as is easily seen by inspection of  standard dihedral group actions on $S^2$. We will use localization to study these actions more closely in Section~\ref{nonabelian}.

With these considerations in mind, let $G= C_p\times C_p=\lan h\ran\times\lan k\ran$, with $H=\lan h\ran$, $K=\lan h\ran$. 
Recall that  $M$ is a four-manifold with $H_1(M;\bz)$ torsion-free,  and  the action of $G$ is HTLL. Assume $\chi(M)\ne 0$.

By Lemma ~\ref{CollapsingLemma1}, the BSS for the action of $H$ collapses, so $H_{H}(M;\bz)\cong H^*(H)\otimes H^*(M)$ as $H^*(H)$-modules. The $G$-isomorphism $H_H^5(M)\cong H^5_H(M^h)$ then gives us a short exact sequence of $G/H$-modules

$$ 0\to H^4( C_p; H^1(M; \bz))\to  H^5_{H}(M^g; \bz)\to H^2( C_p;H^3(M;\bz))\to 0,$$
in which $G/H\cong K $ acts trivially on all terms except possibly the middle.

When $p$ is odd, so that surfaces in $M^h$ are orientable, this becomes 
$$ 0\to \bz_p\otimes H^1(M;\bz)\to \bz_p\otimes H^1(M^h;\bz)\to
\bz_p\otimes H^3(M;\bz)\to 0.$$
(When $p=2$, $H^5_{H}(M^h; \bz)$ contains contributions from the second cohomology of each nonorientable component. )

When $p$ is odd,  it follows from  the classification of integral $ C_p$-representations~\cite[74.3]{CurtisReiner}  that 
$H^1(M^h; \bz)$ splits into a sum of $K$-modules
of trivial, permutation, and cyclotomic types (with no exotic ideal classes; cf.~\cite{SwanGrothendieck}).  

On the other hand, the classification of $\bz_p[K]$-modules  is elementary (cf. ~\cite[2.1]{SikoraTorusZp}): Each splits essentially uniquely as a sum of indecomposibles, and each indecomposible is a cyclic $\bz_p[t]$-module of the form $\bz_p[t]/(t-1)^i$, where $1\le i\le p$. If an indecomposible $\bz_p[K]$-module $M$ arises as $\bz_p\otimes N$, where $N$ is $\bz$-free, then $M$ must have $i=1, p-1$, or $p$. And if $M$ then has both a submodule and a corresponding quotient module with trivial $K$-action, then either $p=2$, or $M$ itself has trivial $K$-action. It follows that

\begin{enumerate}
\item{When $p=2$, $H^1(M^h; \bz_2)$ may contain summands on which $k$ acts by permutation, but}
\item{When $p\ge 3$, the action of $k$ on $H^1(M^h;\bz_p)$ must be trivial.}
\end{enumerate}

\begin{lem}
Let $M$ be a four-manifold whose first homology group is torsion-free, with $\chi(M)\ne 0$, and suppose $M$ admits a homologically trivial, locally linear $ C_p\times C_p$-action, with $p$ an odd prime. Then the singular set of the action consists of isolated points, chains of two-spheres, and isolated tori only.  Moreover, each non-identity group element fixes exactly $b_1(M)$ tori.

\end{lem}
\begin{proof}
Choose any pair of generators $h$, $k$ for the group. The hypotheses guarantee that the BSS for the $\lan h\ran$ action collapses, so $M^h$ consists of points and orientable surfaces. Since the $\lan k\ran$-action on $H^1(M^h)$ must be trivial, the Lefschetz Fixed-point theorem rules out the existence of components of genus $>1$. 
If $T$ is a torus component of $M^h$, and $F$ is a component of $M^k$ meeting $T$, then the $k$-action on $T$ fixes a point, and hence must fix the whole torus. This is impossible, since $G$ cannot act faithfully on a two-dimensional slice disk. So $\lan k\ran$ acts freely on $T$. The torus count then follows from the fact that $b_1(M^h)=2b_1(M)$.

\end{proof}

\begin{thm}\label{OddNoRankTwo}
Let $G= C_p\times C_p$, where $p$ is an odd prime, and let $M$ be a closed topological four-manifold such that $H_1(M;\bz)$ is  torsion-free, with $\chi(M)\ne 0$. If $G$ acts HTLL on $M$, then  $b_1(M)=0$. 
\end{thm}

\begin{proof}
Assume first that $M^G\ne\emptyset$; then it consists of a finite set of points, and by Corollary~\ref{CollapsingLemma2}, the BSS collapses with integer coefficients. Localizing with respect to the set $S\subset H^*(G)$ generated by $\alpha_2, \beta_2\in H^2(G)$ yields an isomorphism $S^{-1}(H^*(G)\otimes H^*(M))\cong S^{-1}(H^*(G)\otimes H^*(M^G) )$. But $S^{-1}H^*(G)$ is a $\bz_2$-graded vector space of dimension $1$ in both even and odd dimensions. It follows that $\sum_{i=0}^4b_i(M)=|M^G|$. It also follows from a simple analysis of the singular set that $\lan k\ran$ acts trivially on $H^*(M^h)$, so by the Lefschetz theorem, $\chi(M) =\chi(M^h)=\chi(M^G) = | \fix(G) |$. Hence $\sum_{i=0}^4b_i(M)=\chi(M)$, so the odd Betti numbers of $M$ must vanish. 

Henceforth we assume (for contradiction) that $M^G=\emptyset$, and that $b_1(M)>0$. Every nontrivial cyclic subgroup $\lan g\ran$ fixes at least one torus. If $g$ were also to fix more than one copy of $S^2$, then by Corollary~\ref{InvariantSurfaces}, all surface components of $M^g$ represent different elements of $H_2(M)$, so each must be $G$-invariant. If $g$ fixes a single $S^2$, then it is obviously $G$-invariant. But a $G$-invariant two-sphere contains a fixed point. The only remaining possibility is that for every $g\in G$, $M^g$ consists of $\chi(M)$ isolated points, freely permuted by $G/\lan g\ran$, and $b_1(M)$ tori, each equipped with a free $G/\lan g\ran$-action. Since $G/\lan g\ran$ freely permutes the isolated points, $p\ |\ \chi(M)$. It follows that $b_2(M)\ge p\ge 3$.

Recall that each cohomology class in $H^2(G; \bz)$ restricts nontrivially to the cohomology of some cyclic subgroup $ C_p\subset G$ -- although the generator $\mu\in H^3(G;\bz)$ does not. Consider the term $E_2^{0,2}$ in the spectral sequence for the $G$-action. Since
$d_2^{0,2}:
H^0(G; H^2(M;\bz))\to H^2(G; H^1(M;\bz)\cong H^2(G)\otimes H^1(M)$, any nontriviality of this differential would be detected by a cyclic subgroup. Hence by Proposition~\ref{CollapsingLemma1}, $d_2^{0,2}=0$.  The differential $d_3^{0,2}$ could be nonzero, but its target $H^3(G; H^0(M))$ is one dimensional. Hence $[E_2^{0,2}:E_{\infty}^{0,2}]\le p$. Since $b_2(M)\ge 3$, there must exist generators $x,y\in H^2(M)$ which survive to $E_{\infty}$ and such that $xy$ generates $H^4(M;\bz)$. Hence $E_2^{0,4}=E_{\infty}^{0,4}$. 

By the multiplicative structure of the spectral sequence, it follows as in  ~\cite[5.4.1]{AlldayPuppeBook} that the $E_{\infty}^{*,4}$ row is a free $H^*(G)$-module , and that the filtration of $H^*_G(M;\bz)$ associated with the spectral sequence yields an exact sequence 
$$0\to\mathcal{F}_3 H^*_G(M;\bz)\to H^*_G(M;\bz)\to E_{\infty}^{*,4}\to 0.$$  Localizing with respect to $S=\lan \alpha, \beta\ran$, it follows that $\fix(G)\ne\emptyset$, a contradiction.

\end{proof}

Unsurprisingly, the case $p=2$ requires special treatment.

\begin{prop} Let $G= C_2\times C_2$ act HTLL on $M$, where $H^*(M;\bz)$ is torsion-free. If $\chi(M)\ne 0$ and $b_2(M)\ne 0,2$, then the BSS collapses with $\bz$ coefficients.
\end{prop}

\begin{proof}
Since $H^2(G;\bz)$ is detected by cyclic subgroups, $d_2$ vanishes. 
For $d_3$, we adapt a calculation from~\cite[6.1]{Actions}:
Let $u\in H^2(M)$. If $u^2\ne 0 \pmod 2$, then $0=d_3(u^3)=3d_3(u)\cdot u^2$, so $d_3(u)=0$. 
On the other hand, if $u^2=0\pmod{2}$, then $b_2(M)\ne 1$, and since we assume $b_2(M)\ne 2$, there exists another generator $v\in H^2(M)$, linearly independent of $u$,  so that $uv=0$. Then $0=d_3(uv)=d_3(u)\cdot v + u\cdot d_3(v).$ But $E_3$ is a free $H^*(G)$-module, so by linear independence,  $d_3(u)=d_3(v)=0$. Hence $d_3$ vanishes on $H^2(M)$. Since $H^4(M)$ is generated by products of two-dimensional classes, $d_3$ vanishes on $H^4$, as well. Finally, let $u$ be a generator of $H^3(M)$. There exists $v\in H^1(M)$ such that $uv$ generates $H^4(M)$, so $d_3(uv)=0$. But $d_3(u)\cdot (v) + u\cdot d_3(v) =0$, and since $d_3(v)=0$, $d_3(u)=0$, as well.

Finally, $d_4$ and $d_5$ vanish on $H^1(M)$ and $H^2(M)$ for dimension reasons, and then on the rest of $H^*(M)$ by Poincar\'e duality. 
\end{proof}

\begin{cor}In the above situation, Localization shows that $\rk H^*(M^G;\bz)=\rk H^*(M;\bz)$.
\end{cor}

Two examples illustrate the additional complications which may arise in the case $p=2$: 

1. When $b_2(M)=2$, $M^{C_2}$ may have two components representing the same homology class, and another involution may interchange them. Indeed, $C_2\times C_2$ acts on $S^2$, and hence on the product of $S^2$ with any surface.

2. 
Let $h=\left(\begin{smallmatrix} 1 & 0& 0&0\\0&1&0&0\\0&0&-1&0\\0&0&0&-1\end{smallmatrix}\right), k=\left(\begin{smallmatrix} 1 & 0& 0&0\\0&-1&0&0\\0&0&-1&0\\0&0&0&1\end{smallmatrix}\right), $
restricted to $S^3\subset\br^4$,
and let $M$ be the mapping cylinder $S^3\times I/\{(x,0)\sim (hk(x), 1)\}$. 
Since $hk$ is isotopic to the identity, $M\cong S^3\times S^1$. Since each of $h$ and $k$ commutes with $hk$, there is an induced $\bz_2\times\bz_2$ action on $M$ with $M^h\cong M^k$ each homeomorphic to a Klein bottle, while $T^2=M^{hk}$ is a torus. The global fixed set is a pair of circles.

Let us say that an action of an elementary 2-group $G$ on a spin 4-manifold is \emph{weakly spin-preserving} (WSP) if at least one involution $g\in G$ preserves a spin structure on $M$. (Smooth vs LL?) It follows from~\cite{EdmondsOrientability, Aspects} that $\fix{g}$ is orientable.

Example 2 is WSP. Indeed, $hk$ preserves both spin structures on $M$, while each of $h$ and $k$ interchanges them. Each of $h$ and $k$ restricts to an orientation-reversing reflection on $T^2$. Such behavior would be ruled out by the assumption of homological triviality if $[T^2]$ were nontrivial in $H_2(M;\bz)$. 

When $b_2(M)\ne 0,2$, each component $F$ of $\fix(h)$ is $G-$invariant and nontrivial in $H_2(M, \bz)$. Hence any $k\in G$ will preserve the orientation of $F$ and act by rotation in the neighborhood of a fixed point, so $\fix(h)\cap \fix(k)$ will be a discrete set of points. 

Another situation in which orientability is ensured occurs when $G$ acts by symplectic symmetries (cf.~\cite{ChenKwasikSymplecticSymmetries}) on a symplectic four-manifold. Components of the fixed-point set are symplectic submanifolds, and are therefore oriented.

\begin{thm}

Let $G= C_2\times C_2$,  and let $M$ be a closed topological four-manifold such that $H^*(M;\bz)$ is  torsion-free, with $\chi(M)\ne 0$. and $b_2(M)\ne 0,2$.  If $G$ acts HTLL  and either
\begin{enumerate}
	\item{$M$ is spin, and the $G$-action is WSP, or}
	\item{$M$ is symplectic, and $G$ acts by symplectic symmetries,}
	\end{enumerate}
	 then  $b_1(M)=0$. 
\end{thm}
\goodbreak

\begin{proof}\label{Z2xZ2Theorem}
Let $G=\lan h, k\ran$, where $h$ preserves a spin structure or a symplectic 2-form.  In the proof of Theorem~\ref{OddNoRankTwo} for odd $p$, we knew the $G$-action on  $H^1(M^h)$ was trivial, and in the presence of a fixed point, the result followed immediately from the Lefschetz theorem. Here, we have not ruled out the possibility that $h$ might act by permutations on part of $H^1(M^h)$.

It is a consequence of localization (cf. ~\cite[III.4.16]{tomDieck}) that for any $ C_p$-space $X$ such that $\dim_{\bz_p}\bigoplus H^k(X; \bz_p)$ is finite, we have 
$$\dim_{\bz_p}\bigoplus H^k(X^{ C_p}; \bz_p)\le \dim_{\bz_p}\bigoplus H^k(X; \bz_p),$$
with equality if and only only if $ C_p$ acts trivially on $H^*(X;\bz_p)$ and the BSS collapses. 

Separate $M^h$ into its $\lan k\ran$-orbits of path-components, say $X_1, \ldots, X_n$. Each is either an oriented surface, a point, or a pair of points. Then \begin{align*}
|\ M^G\ | &=\sum_{i=1}^n |\ ( X_i)^k\ |\\
&\le \sum_{i=1}^n (\dim_{\bz_p}\bigoplus H^j(X_i;\bz_p))\\
&=\dim_{\bz_p}\bigoplus H^i(M^h;\bz_p)\\
&= \sum b_j(M) \\
&=|\ M^G\ |.
\end{align*}

It follows that for each individual $X_i$, $|(X_i)^k|=\sum(b_j(X_i))$.
By the Lefschetz Theorem, this is only possible if $\trace(k|_{H^1(X_i)})=-b_1(X_i)$. But the trace of a permutation representation is zero. Hence $b_i(X_i)=0$ for all $X_i\subset M^h$. So $b_1(M^h)=0$, and finally, $b_1(M)=0$.
\end{proof}

Without any assumption to guarantee orientability of the fixed-point set, a corank two subgroup still acts on a collection of points and circles, with $\sum b_i(M^G) = \sum b_i (M)$, so $G$ can still be no larger than $(C_2)^3$. Finer analysis might make it possible to remove the assumption, but the presence of nonorientable surfaces in the singular set does complicate its combinatorial structure.

\section{Nonabelian groups}\label{nonabelian}

Let $p$ be prime. Recall that the $p$-\emph{rank} of a finite group $G$ is the largest $n$ such  $( C_p)^n\subset G$. It follows from the results of the last section that if $M$ is a closed,  oriented four-manifold with $H_1(M)$ nontrivial and torsion-free, $\chi(M)\ne 0$, $b_2(M)\ne 0,2$, then any finite group which acts HTLL on $M$ must have $p$-rank $\le 1$ for any odd prime $p$.   The corresponding assertion for $p=2$ also assumes the action to be weakly spin-preserving.  Without this assumption, we know $G$ has 2-rank $\le 3$.

In ~\cite{MM2}, we analyzed minimal nonabelian groups and saw in particular that if $G$  is a nonabelian rank one finite group such that every proper subgroup of $G$ is abelian, then $G$ is either a metacyclic group of the form $ C_p\rtimes C_{q^n}$, where $p$ and $q$ are prime, and $ C_{q^n}$ acts on $ C_p$ via an order $q$ group automorphism, or $Q_8$, the order 8 quaternion group. 

If $G$ is a minimal nonabelian group with $2$-rank $\le 3$, and $p$-rank $\le 1$ for odd $p$, there are additional possibilities: 2-groups of the form $G_1(m,n,2) = \lan a,b |\ a^{2^m}=b^{2^n}=1, [a,b]=a^{2^m-1}\ran,$ for $m\ge 2$; $G_2(m,n,2) = \lan a,b,c |\ a^{2^m}=b^{2^m}=c^2=1, [a,b]=c, [c,a] = [c,b]=1\ran$. It turns out that $G_1(2,1,2)\cong G_2(1,1,2)\cong D_4$; each of the others in these families contains a subgroup isomorphic to $C_4\times C_2$. Finally, there are $(C_2\times C_2)\rtimes C_3 = A_4$, and $(C_2\times C_2\times C_2)\rtimes C_7$. \emph{Ad hoc} arguments to rule out actions of these groups are sometimes possible. 

Our remaining arguments rely on localization to yield information about the action of a group on the fixed set of a normal subgroup.  As applications in the literature seem to focus mainly on abelian groups and connected Lie groups,  we begin with a simple example for motivation. Let $p$ be an odd prime, and consider a standard linear action of the dihedral group $D_p=\lan a,b\ |\ a^p=b^2=1, bab^{-1}=a^{-1}\ran$ on $S^2$. Each cyclic subgroup fixes two points; each involution interchanges the points in $(S^2)^a$, and the action of $\lan a\ran =C_p$ permutes the fixed sets of the various $2$-subgroups.

The integral cohomology of $D_p$ is periodic of period 4, vanishing in odd dimensions, cyclic of order 2 in dimensions equal to 2 mod 4, and cyclic of order $2p$ in dimensions divisible by 4. (This is well-known, but also follows from our general calculation for metacyclic groups below). With this in mind it is easy to see that the BSS for the $D_p$ action on $S^2$ collapses, so $H^*_{D_p}(S^2;\bz)$ is a free $H^*(D_p)$ module on generators of degree zero and two.

In the case of our $D_p$ action, the periodicity generator $u$ is detected by restriction to both $ C_p$ and $ C_2$, so if $S$ is the multiplicative set generated by $u$, then $\Sigma_S$ is the entire singular set. Hence localization tells us nothing we did not already know by simply considering the (non-localized) equivariant cohomology groups in degrees $d>2$. However, if $S$ is generated instead by $2u$, then $\Sigma_S =\fix( C_p)$, and $S^{-1}H^*(D_p)$ is a $\bz_4$-graded module with $ \bz_p$ in degree four, and zeroes otherwise. It follows that as a $S^{-1}H^*(D_p)$-module, 
$S^{-1}H_{D_p}^*(\fix( C_p))$ has one $\bz_p$ in each even degree. 

We know (geometrically) that $\fix( C_p)$ is a $0$-dimensional manifold, equipped with an action of $D_p/ C_p\cong C_2$. It follows that $H^0(\fix( C_p))$ is a $\bz$-free $\bz[ C_2]$-module, and hence a sum of modules of trivial, permutation, and cyclotomic types. But one checks that $H^*(D_p;\bz[D_p/ C_p])$ has generators in every even degree, and $H^*(D_p;\bz[-1])$ has generators in degrees $d\cong 2\pmod{4}$. Localization therefore detects the fact that that the action of $D_p/ C_p$ on $\fix( C_p)$  is nontrivial: either entirely by permutations, or with equally many trivial and cyclotomic components. It is simple to rule out the latter case with other considerations, but even without any, note the consequence that the trace of the $ C_2$-action on $H^0(\fix( C_p))$ is zero. 

Now consider the case of a metacyclic group $G=\lan a,b\ |\ a^p=1=b^{q^n}, bab^{-1}=b^r\ran$, where $p$ and $q$ are prime, and $r^q\equiv 1\pmod p$. 
The Lyndon-Hochschild-Serre spectral sequence of the group extension
$1\to C_p\to G\to C_{q^n}\to 1$ has $H^i( C_{q^n}; H^j( C_p;\bz))\Rightarrow H^{i+j}(G;\bz)$. Since $p$ and $q$ are relatively prime, $H^i(C_{q^n}; H^j( C_p;\bz))=0$ whenever $i$ and $j$ are both positive. And in general, $H^0( C_{q^n}, X)\cong X^{C_{q^n}}$ for any $ C_q$-module $X$. The action of $C_{q^n}$ on $H^*(C_p;\bz)$ is determined by its action on the generator $t$ of $H^2(C_p)$, and $b\cdot t=rt$. Hence \begin{equation*}
H^0( C_{q^n}; H^j( C_p;\bz))=\begin{cases}
\bz &\text{ if $j=0$,}\\
\bz_p &\text{if $j>0$ and $2q\ |\ j$},\\
0 &\text{otherwise.}
\end{cases}
\end{equation*}
Of course, $H^i( C_{q^n};H^0( C_p;\bz))=\bz_{q^n}$ in even dimensions, and zero in odd. The spectral sequence collapses, and so we see that
 \begin{equation*}
H^i(G;\bz)=\begin{cases}
\bz &\text{if $i=0$,}\\
\bz_{q^n} &\text{if $i$ is even, but $q\not{|}\ i$,}\\
\bz_{pq^n} &\text{if $2q\ |\  i$,}\\
0 &\text{otherwise},
\end{cases}
\end{equation*}
and that restriction to the subgroup $\bz_{q^n}$ is a cohomology isomorphism in dimensions below $2q$.

\begin{lem} \label{MetacyclicCollapsingLemma}Let $M$ be a closed, oriented four-manifold with $H_1(M)$ torsion-free, and $\chi(M)\ne 0$. If a nonabelian group of the form $G= C_p\rtimes C_{q}$ (with $p$ and $q$ prime) acts HTLL on $M$, then the $\bz$-coefficient BSS for the action collapses. 
\end{lem}

\begin{proof} Pick a representative subgroup $ C_q\subset G$, If $q>2$, then the target of every differential $d_r^{0,j}$
lies within the region where restriction to the cohomology of $ C_q$ is an isomorphism, and factoring through the spectral sequence for $H_{ C_q}^*(M)$ shows that the given sequence collapses. 

If $q=2$, then the differential $d_4^{0,3}:H^0(G;H^3(M;\bz))\to H^4(G; H^0(M;\bz)\cong  \bz_{2p}$ needs special consideration. But restriction to $C_2$ takes care of potential $2$-torsion, and restriction to $ C_p$ rules out $p$-torsion.
\end{proof}

\begin{thm}Let $M$ be a closed, oriented four-manifold with $H_1(M)$ nontrivial and torsion-free,  $\chi(M)\ne 0$, and $b_2(M)\ne 0,2$. The group $G= C_p\rtimes C_{q}$ (with $p$ and $q$ prime) admits no effective,  HTLL action on $M$. \end{thm}

\begin{proof}Suppose such an action exists. By Lemma~\ref{MetacyclicCollapsingLemma}, $H_{G}^*(M;\bz)\cong H^*(G)\otimes H^*(M;\bz)$. Let $u\in H^{2q}(G;\bz)$ be a periodicity generator, and consider $S=\lan qu\ran\subset H^*(G;\bz)$. In this case, $\Sigma_S=M^{C_p}$, and $S^{-1}H^*(G)$ is a $\bz_{2q}$-graded module which is  isomorphic to $\bz_p$ in dimensions divisible by $2q$, and trivial otherwise. It follows that 
\begin{equation*}
S^{-1}H_G^{\text{odd}}(M^{C_p})\cong
\begin{cases}(\bz_p)^{b_1(M)} &\text{ in degrees $d\equiv1$ and $d\equiv 3\pmod {2q}$,}\\
0&\text{ otherwise.}
\end{cases}
\end{equation*}
However,  $S^{-1}H_G^{\text{odd}}(M^{C_p})\cong S^{-1}H^{2q}(G; H^1(M^{C_p})$.
Now, $H^1(M^{C_p})$ is a $\bz$-free $\bz[ C_q]$-module, and hence equivalent to a sum of modules of cyclotomic, trivial, and permutation types. By Shapiro's lemma, $H^*(G; \bz[ C_q])\cong H^*( C_p;\bz)$, and in the cases $d=2qk$, the restriction $H^{d}(G;\bz)\to H^{d}( C_p;\bz)$ is onto. It follows from the coefficient short exact sequence $0\to \bz[\lambda]\to\bz[ C_q]\to \bz\to 0$  that
\begin{equation*}
H^i(G;\bz[\lambda])=
\begin{cases}
\bz_p &\text {if $i$ is even and $2q\not{|}\ i$,}\\
\bz_q &\text{if $i=2qk+1$,}\\
0&\text{otherwise.}
\end{cases}
\end{equation*}
so each cyclotomic summand in $H^1(M^{C_p})$ contributes $\bz_p$ to $S^{-1}H_G^*(M^{C_p})$ in each odd degree $d$ except $d\equiv 1\pmod{2q}$. Each free summand makes a contribution in every odd degree, and each trivial summand makes a contribution when $d\equiv 1\pmod{2q}$.

If $q> 2$, then $S^{-1}H^5_G(M^{C_p})=0$, which rules out the presence of cyclotomic and permutation summands in $H^1(M^{C_p})$. On the other hand, trivial summands cannot account for the nontriviality of $S^{-1}H^3_G(M^{C_p})$. It follows that no actions can exist. 

If $q=2$, a similar analysis shows that the numbers of cyclotomic and trivial summands in $H_1(M^{C_p})$ are equal, and it follows that the trace of the $ C_2$ action on $H_1(M^{C_p})$ is zero. Since $b_0(M^{C_p})+b_2(M^{C_p})=b_2(M)+2$, and $b_2(M)\ne 0,2$, $M^{C_p}$ has either a single two-dimensional component, or its distinct components represent different elements of $H_2(M)$.  So the two-dimensional components of $M^{C_p}$ are individually $G$-invariant. It follows that the Lefschetz number of the $C_2$-action on $M^{C_p}$ is positive, and that a fixed point $x$ for the entire group action exists on some two-dimensional component of $M^{C_p}$. However, 
any $2$-dimensional component $F$ of $M^{C_p}$ forms a proper subset of $M^{C_p}$, and hence each such component represents a nontrivial class in $H_2(M; \bz)$. Consideration of the local representation of the dihedral group $G$ on $T_x(M)$ shows that the $ C_2$-action must reverse orientation on the $2$-plane fixed by $ C_p$, and hence send $[F]$ to $-[F]$, contradicting homological triviality.

\end{proof}

A refinement of this argument rules out the larger metacyclic groups:

\begin{thm}Let $M$ be a closed, oriented four-manifold with $H_1(M)$ nontrivial and torsion-free, $\chi(M)\ne 0$, and $b_2(M)\ne 0,2$. The group $G=C_p\rtimes C_{q^n}$ (with $p$ and $q$ prime, and $n>1$) admits no effective,  HTLL action on $M$. \end{thm}

\begin{proof}

Let $G=\lan a,b\ |\ a^p=1=b^{q^n}, bab^{-1}=b^r\ran$, where $r^q\cong 1\pmod p$, and begin by considering the index $q$ cyclic subgroup $H=\lan ab^q\ran$. The argument of the previous section shows that $\lan b^q\ran$ acts on the spectral sequence for the $ C_p$ action, and hence we have a short exact sequence of $ C_{q^{n-1}}$-modules
$$ 0\to H^4( C_p; H^1(M))\to \bz_p\otimes H^1(M^a)\to H^2( C_p;H^3(M))\to 0,$$
where $ C_{q^{n-1}}$ acts trivially on the outside terms. By Maschke's theorem,  $ \bz_p\otimes H^1(M^a)$ must have trivial $ C_{q^{n-1}}$-action. Hence $H$ acts trivially on $H^1(M^{C_p};\bz_p)$, and the $ C_{q^n}$ action on $\bz_p\otimes H^1(M^{C_p})$ factors through a $G/H\cong  C_q$-action. 

Now consider the spectral sequence of the $G$-action on $M$. It is not evident \emph{a priori\ } that it collapses, but factoring through the BSS for the the $C_p$-action shows that that all $p$-torsion survives, and it follows again that with $S=\lan qu\ran$, we still have 
\begin{equation*}
S^{-1}H_G^{\text{odd}}(M^{C_p})\cong
\begin{cases}(\bz_p)^{b_1(M)} &\text{ in degrees $d\equiv1$ and $d\equiv 3\pmod {2q}$,}\\
0&\text{ otherwise.}
\end{cases}
\end{equation*}
The rest of the argument proceeds exactly as before. 
\end{proof}

Finally, we turn our attention to the quaternion group $$Q_8=\lan h, x\ |\ h^4=x^4=1, xhx^{-1}=h^{-1}, h^2=x^2\ran.$$ $Q_8$ has three cyclic subgroups of order 4: $\lan h\ran, \lan x\ran,$ and $\lan xh\ran$, all of which intersect in the central $\lan h^2\ran\cong  C_2$. 

\begin{thm} Let $M$ be a closed, oriented four-manifold with $H_1(M)$ nontrivial and torsion-free, and $\chi(M)\ne 0$. The group $Q_8$ admits no effective,  HTLL action on $M$. \end{thm}

\begin{proof}
Suppose $Q_8$ were to act in the stated manner. Denote the three copies of $ C_4\subset Q_8$ by $K_1, K_2, K_3$, and 
consider the action of a subgroup $K_i$. The spectral sequence for the $K_i$-action (with $\bz_2$ coefficients) may not collapse, but by Proposition~\ref{CollapsingLemma1}, the only potentially nonzero differential is $d_2$ from row 2 to row 1. In particular, row 3 consists of permanent cocycles.

Now, examination of the L-H-S spectral sequence (for example) shows that the restriction map $H^i( C_4;\bz_2)\to H^i( C_2;\bz_2)$ is trivial for every $i>0$. Hence if $u$ generates $H^2( C_4;\bz_2)$, and $S=\lan u\ran$, then $\Sigma_S$ consists only of points fixed by the entire group $ C_4$. Hence, by the Localization Theorem, $b_1(M^{K_i})\ge b_1(M)$.  

But each $K_i\subset Q_8$ has the center $ C_2$ as a subgroup, and clearly, $M^{K_i}\subset M^{C_2}$. Let $F$ be the set of two-dimensional components of $M^{C_2})$, and for each $i$, let $F_i$ be the corresponding set of two-dimensional components of $M^{K_i}$. Let $V=H_1(F;\bz_2)$, and $V_i=H_1(F_i;\bz_2)\subset V$. By a dimension count, it follows that some pair $V_i, V_j$ must have nontrivial intersection, and hence that $F_i$ and $F_j$ have at least one two-dimensional component in common. But such a component is necessarily fixed by all of $Q_8$. This is impossible, as $Q_8$ cannot act faithfully on a two-dimensional slice disk. 
\end{proof}

As we have eliminated groups of rank $\ge 2$ and all nonabelian groups, we have established the main theorem.

\bibliographystyle{abbrv}
\bibliography{mybiblio}

\begin{thebibliography}{10}

\bibitem{AlldayPuppeBook}
C.~Allday and V.~Puppe.
\newblock {\em Cohomological methods in transformation groups}, volume~32 of
  {\em Cambridge Studies in Advanced Mathematics}.
\newblock Cambridge University Press, Cambridge, 1993.

\bibitem{Assadi}
A.~H. Assadi.
\newblock Integral representations of finite transformation groups {III}:
  Simply-connected four-manifolds.
\newblock {\em Journal of Pure and Applied Algebra}, 65(3):209--233, 1990.

\bibitem{Borel}
A.~Borel.
\newblock {\em Seminar on Transformation Groups}, volume~46 of {\em Annals of
  Mathematics Studies}.
\newblock Princeton University Press, 1960.

\bibitem{Bredon}
G.~E. Bredon.
\newblock {\em Introduction to Compact Transformation Groups}.
\newblock Academic Press, New York, 1972.

\bibitem{BredonSheafTheory}
G.~E. Bredon.
\newblock {\em Sheaf theory}, volume 170 of {\em Graduate Texts in
  Mathematics}.
\newblock Springer-Verlag, New York, second edition, 1997.

\bibitem{ChenKwasikSymplecticSymmetries}
W.~Chen and S.~Kwasik.
\newblock Symplectic symmetries of 4-manifolds.
\newblock {\em Topology}, 46(2):103--128, 2007.

\bibitem{CurtisReiner}
C.~W. Curtis and I.~Reiner.
\newblock {\em Representation theory of finite groups and associative
  algebras}.
\newblock Pure and Applied Mathematics, Vol. XI. Interscience Publishers, a
  division of John Wiley \& Sons, New York-London, 1962.

\bibitem{EdmondsOrientability}
A.~L. Edmonds.
\newblock Orientability of fixed point sets.
\newblock {\em Proc. Amer. Math. Soc.}, 82(1):120--124, 1981.

\bibitem{Aspects}
A.~L. Edmonds.
\newblock Aspects of group actions on four-manifolds.
\newblock {\em Topology and its Applications}, 31(2):109--124, 1989.

\bibitem{Actions}
A.~L. Edmonds.
\newblock Homologically trivial group actions on 4-manifolds.
\newblock Preprint; electronically posted at
  http://front.math.ucdavis.edu/math.GT/9809055, 1997.

\bibitem{Hsiang}
W.~Y. Hsiang.
\newblock {\em Cohomology Theory of Topological Transformation Groups},
  volume~85 of {\em Ergebnisse der Mathematik und ihrer Grenzgebiete}.
\newblock Springer-Verlag, New York, 1975.

\bibitem{Lewis}
G.~Lewis.
\newblock The integral cohomology rings of groups of order $p^3$.
\newblock {\em Transactions of the American Mathematical Society},
  132:501--529, 1968.

\bibitem{MM2}
M.~P. McCooey.
\newblock Symmetry groups of four-manifolds.
\newblock {\em Topology}, 41(4):835--851, 2002.

\bibitem{SikoraTorusZp}
A.~S. Sikora.
\newblock Torus and {${\Bbb Z}/p$} actions on manifolds.
\newblock {\em Topology}, 43(3):725--748, 2004.

\bibitem{SwanGrothendieck}
R.~G. Swan.
\newblock The {G}rothendieck ring of a finite group.
\newblock {\em Topology}, 2:85--110, 1963.

\bibitem{tomDieck}
T.~tom Dieck.
\newblock {\em Transformation Groups}.
\newblock Walter de Gruyter \& Co., Berlin, 1987.

\end{thebibliography}

\end{document}